\documentclass[12pt,a4paper]{amsart}
\usepackage{amssymb, enumerate, mdwlist, amsthm, amsmath, bm, graphicx}
\usepackage{graphicx,epsfig}
\usepackage[dvipdfmx]{xcolor,hyperref}
\hypersetup{
    colorlinks=true,
    citecolor=blue,
    linkcolor=red,
    urlcolor=orange,
}
\newtheorem{thm}{Theorem}
\newtheorem{prop}{Proposition}
\newtheorem{conj}{Conjecture}
\newtheorem{lemma}{Lemma}
\newtheorem{claim}{Claim}
\def\Z{\mathbb Z}
\def\AA{\mathcal A}
\def\FF{\mathcal F}
\def\GG{\mathcal G}
\def\HH{\mathcal H}
\def\WW{\mathcal W}

\title{The maximum measure of 3-wise $t$-intersecting families}
\author{Norihide Tokushige}
\address{University of the Ryukyus\\College of Education\\
1 Senbaru Nishihara\\Okinawa\\903-0213 (JAPAN)\\
ORCID: 0000-0002-9487-7545}
\email{hide@edu.u-ryukyu.ac.jp}
\date{\today}

\begin{document}
\maketitle
\begin{abstract}
Let $\GG$ be a family of subsets of an $n$-element set. 
The family $\GG$ is called $3$-wise $t$-intersecting if the intersection of 
any three subsets in $\GG$ is of size at least $t$.
For a real number $p\in(0,1)$ we define the measure of the family 
by the sum of $p^{|G|}(1-p)^{n-|G|}$ over all $G\in\GG$. 
We prove that if $t\geq 15$ and $p\leq 2/(\sqrt{4t+9}-1)$ then $p^t$ is 
the maximum measure of $3$-wise $t$-intersecting families, 
and the bound for $p$ is sharp. 
We also present the corresponding stability result for shifted families.
\end{abstract}

\section{Introduction}
Let $n\geq t\geq 1$ and $r\geq 2$ be integers.
Let $[n]=\{1,2,\ldots,n\}$ and let $2^{[n]}$ denote the power set of $[n]$.
We say that a family $\GG$ of subsets is
$r$-wise $t$-intersecting if $|G_1\cap\cdots\cap G_r|\geq t$ for all
$G_1,\ldots,G_r\in\GG$. If $t=1$ then we omit $t$ and say an $r$-wise 
intersecting family to mean an $r$-wise $1$-intersecting family.

Let $0<p<1$ be a real number and let $q=1-p$. 
For $\GG\subset 2^{[n]}$ we define its measure (or $p$-biased measure) 
$\mu_p(\GG)$ by
\[
 \mu_p(\GG):=\sum_{G\in\GG}p^{|G|} q^{n-|G|}.
\]
There are two basic problems concerning how large an $r$-cross 
$t$-intersecting family can be. The first one asks the maximum size of
$k$-uniform families, that is,
\[
\max\{|\FF|:\FF\subset\binom{[n]}k\text { is $r$-wise $t$-intersecting}\}.
\]
The second one asks the maximum $p$-biased measure, that is,
\[
\max\{\mu_p(\FF):\FF\subset 2^{[n]}\text { is $r$-wise $t$-intersecting}\}.
\]
These problems are closely related, in particular, when $p\approx\frac kn$.
More generally, some results in the $k$-uniform setting transfer to
the corresponding result in the $p$-biased setting, and vice versa,
see e.g., \cite{EKL18, EKL19, Friedgut, Tuvsw}.

With respect to these two problems, it seems natural to 
expect that the maximum is attained by a family of subsets
(in $\binom{[n]}k$ or $2^{[n]}$) containing at least $t+(r-1)s$ points 
out of some fixed set of $t+rs$ points.
These problems have a long history going back to Erd\H{o}s--Ko--Rado
\cite{EKR}, and have been extensively studied. However, the complete 
solutions are known only for the two cases in both settings: 
$r=2$ and $t\geq 1$ (\cite{AK1, AK2, BE} for $k$-uniform setting, and
\cite{AK-p, DS, Filmus, Friedgut, Tuvsw} for $p$-biased setting), 
and $r\geq 2$ and $t=1$ (\cite{F76, F87, Gronau, MV} for 
$k$-uniform setting, and \cite{FGL, FT2003, T2021} for $p$-biased setting.)
Only a few results are known for the other cases $r\geq 3$ and $t\geq 2$.
The cases in the $p$-biased setting where 
$p=\frac12$ (\cite{Frankl1991,Frankl2019}) and
$p$ is close to $\frac12$ (\cite{T2007b}), and the cases where
$r\geq r_0(p)$ and $t\leq t_0(p,r)$ in both settings (\cite{T2010}).
In this paper we deal with the case $r=3$ and $t\geq 15$ in the $p$-biased 
setting.

For $s=0,1,\ldots\lfloor\frac{n-t}3\rfloor$, 
we define a 3-wise $t$-intersecting family 
$\FF_s^t=\FF_s^t(n)\subset 2^{[n]}$ by
\[
 \FF_s^t:=\{F\in 2^{[n]}:|F\cap[t+3s]|\geq t+2s\}.
\]
Then $\mu_p(\FF_s^t)=\sum_{j=t+2s}^{t+3s} \binom{t+3s}j p^{j}q^{t+3s-j}=\sum_{i=0}^s \binom{t+3s}i p^{t+3s-i}q^i$, 
and a direct computation shows that 
$\mu_p(\FF_0^t)\geq\mu_p(\FF_1^t)$ if and only if $p\leq p_0$, where
\begin{align}\label{def:p0}
 p_0=p_0(t):=\frac 2{\sqrt{4t+9}-1}. 
\end{align}
We say that two families $\FF,\FF'\subset 2^{[n]}$ are isomorphic, denoted by 
$\FF\cong\FF'$, if there is a permutation $\pi$ on $[n]$ such that
$\FF'=\{\{\pi(x):x\in F\}:F\in\FF\}$. Our main result is the following.

\begin{thm}\label{thm1}
Let $t\geq 15$ and $0<p\leq p_0$.
Let $\GG\subset 2^{[n]}$ be a $3$-wise $t$-intersecting family.
Then $\mu_p(\GG)\leq p^t$. Moreover equality holds if and only if 
$\GG\cong\FF_0^t$, or $p=p_0$ and $\GG\cong\FF_1^t$.
\end{thm}
In \cite{T2007} one can find the corresponding result in $k$-uniform
setting for $t\geq 26$.
Both proofs are based on the same method described in the next section, 
but the casewise analysis is different. For the proof of Theorem~\ref{thm1}
we divide $\GG$ into three subfamilies, and estimate the $p$-biased measure
of each subfamily separately. As a result we get a simpler proof with a 
better lower bound on $t$.

We conjecture that the same conclusion as in Theorem~\ref{thm1}
holds for $2\leq t\leq 14$ as well.
The situation for the case $t=1$ is different. On one hand
if $s=\lfloor\frac{n-1}3\rfloor$ then 
$\lim_{n\to\infty}\mu_p(\FF^1_s)=1$ for $p>\frac23$. 
On the other hand it is known from \cite{FGL,FT2003,T2021} that every 
3-wise 1-intersecting family $\FF$ satisfies $\mu_p(\FF)\leq p$ for 
$0<p\leq\frac23$, and if $p<\frac23$ then $\mu_p(\FF)=p$ holds only if 
$\GG\cong\FF_0^1$.

We also prove the following stability result, which says that if a
3-wise intersecting family has $p$-biased measure very close to the optimal
value $p^t$, then the family itself is close to the optimal family ($\FF_0^t$
or $\FF_1^t$) in structure. Here we consider two families $\FF$ and $\GG$
to be close in structure if the symmetric difference 
$\FF\triangle\GG:=(\FF\cup\GG)\setminus(\FF\cap\GG)$ has small 
$p$-biased measure.
We say that a family $\FF\subset 2^{[n]}$ is shifted if 
$(F\setminus\{j\})\cup\{i\}\in\FF$ for all
$F\in\FF$ satisfying $\{i,j\}\cap F=\{j\}$ for some $1\leq i<j\leq n$.

\begin{thm}\label{thm2}
Let $t\geq 15$ and $0<p\leq p_0$.
Then there exist constants $\epsilon_0=\epsilon_0(p,t)>0$ and
$C=C(p,t)>0$ satisfying the following statement: for every 
$0<\epsilon<\epsilon_0$ and every shifted $3$-wise $t$-intersecting family
$\GG\subset 2^{[n]}$ with $\mu_p(\GG)=p^t-\epsilon$, it follows that
\[
 \mu_p(\FF\triangle\GG)<C\epsilon,
\]
where $\FF\in\{\FF_0^t,\FF_1^t\}$.
\end{thm}

If $t$ is sufficiently large, then we have the following stronger stability.

\begin{thm}\label{thm4}
There exists $t_0$ such that for every $t\geq t_0$ and every $p<p_0$ 
there exists a constant $\epsilon_0=\epsilon_0(p,t)>0$ satisfying
the following statement: for every $0<\epsilon<\epsilon_0$
and every shifted $3$-wise $t$-intersecting family
$\GG\subset 2^{[n]}$ with $\mu_p(\GG)=p^t-\epsilon$, we have 
\[
 \mu_p(\GG\setminus\FF)<\frac{3\epsilon}{t-3},
\]
where $\FF\in\{\FF_0^t,\FF_1^t\}$.
\end{thm}

Friedgut obtained a similar result for 2-wise $t$-intersecting families
in \cite{Friedgut}, which was one of the earliest results concerning 
stability of intersecting families. His proof is based on a deep result, 
due to Kindler and Safra \cite{KS}, concerning
Boolean functions whose Fourier transforms are concentrated on small sets. 
In the same approach with tools from \cite{FGL}, 
a stability result for 3-wise 1-intersecting families was obtained in
\cite{T2021}. 
Ellis, Keller and Lifshitz \cite{EKL18} proved stability of the
Ahlswede--Khachatrian theorem \cite{AK1} by modifying its proof,
combining combinatorial techniques with a simple isoperimetric result.
All these stability results in \cite{Friedgut, T2021, EKL18} 
are stronger in the sense that they do not need to assume that the family 
is shifted.

Our proof of Theorem~\ref{thm2} is purely combinatorial. 
Applying the same proof method to $k$-uniform families (as in
\cite{T2007}), it seems possible to obtain the results corresponding to 
Theorems~\ref{thm2} and \ref{thm4} in the $k$-uniform setting.
However, it does not seem possible to directly derive such stability 
results in the $k$-uniform setting from Theorems~\ref{thm2} and \ref{thm4}.

Finally we conjecture a much stronger stability for all $t\geq 2$.

\begin{conj}
Let $t\geq 2$ and $0<p\leq p_0$.
Then there exists $\epsilon=\epsilon(p,t)>0$ satisfying the following 
statement: if $\GG\subset 2^{[n]}$ is a $3$-wise $t$-intersecting family with
$\mu_p(\GG)>p^t-\epsilon$, then $\GG\subset\FF$ for some $\FF\subset 2^{[n]}$
such that $\FF\cong\FF_0^t$, or $p=p_0$ and $\FF\cong\FF_1^t$.
\end{conj}

The corresponding result for $2$-wise $t$-intersecting families is not
true. To see this, let 
$\FF':=(\FF_0^t\setminus\{[t]\})\sqcup\{[n]\setminus\{1\}\}\subset 2^{[n]}$. 
Then $\FF'\not\subset\FF_0^t$ (and $\FF'\not\subset\FF_1^t$) but
$\mu_p(\FF')=p^t-p^tq^{n-t}+p^{n-1}q=p^t+o(1)$ as $n\to\infty$.
Related examples are given in the discussion after the following
theorems: Theorem 1.4 in \cite{FLST}, Theorem 1.9 in \cite{EKL19}, 
Theorem 1.10 in \cite{EKL18}, and Theorem 3 in \cite{LST}.
See also Theorem 6 in \cite{T2022} for positive evidence for shifted 
non-trivial 3-wise intersecting families.

\section{Preliminaries}
For the proof of Theorem~\ref{thm1} we use the random walk method introduced
by Frankl in \cite{F87}, where the case $p=\frac12$ is considered.
In this section we gather some basic facts with proofs for convenience. 
See \cite{LST} and the references therein for more about 
recent development of the method.

Let $\GG\subset 2^{[n]}$.
For a subset $G\in\GG$ we define its walk in $\Z^2$.
This is an $n$-step walk starting from the origin, and the $i$th step is up
if $i\in G$ or right if $i\not\in G$. We will identify $G$ and its walk, and
we say that $G$ hits a line to mean that the walk of $G$ hits the line.
(Here line hitting includes line crossing.)
So $\FF_s^t$ is the set of $n$-step walks hitting one of 
$(0,t+3s), (1,t+3s-1),\ldots, (s,t+2s)$, e.g.,
$\FF_1^t$ is the set of $n$-step walks hitting $(0,t+3)$ or $(1,t+2)$.

For $1\leq i<j\leq n$ we define the shifting operation 
$\sigma_{i,j}:2^{[n]}\to 2^{[n]}$ by 
\[
\sigma_{i,j}(\GG):=\{G_{i,j}:G\in\GG\}, 
\]
where
\[
G_{i,j}:=\begin{cases}
(G\setminus\{j\})\cup\{i\}&\text{if }G\cap\{i,j\}=\{j\}
\text{ and }(G\setminus\{j\})\cup\{i\}\not\in\GG,\\
G & \text{otherwise}.	   
\end{cases}
\]
By definition $\mu_p(\GG)=\mu_p(\sigma_{i,j}(\GG))$.
We say that $\GG$ is shifted if $\GG$ is invariant under any shifting 
operation, in other words, if $G\in\GG$ then $G_{i,j}\in\GG$ for all 
$1\leq i<j\leq n$. If $\GG$ is not shifted then 
$\sum_{G\in\GG}\sum_{g\in G}g>\sum_{G'\in\sigma_{i,j}(\GG)}\sum_{g'\in G'}g'$
for some $i,j$, and so starting from $\GG$ we get a shifted $\GG'$ by 
applying shifting operations a finite number of times.
It is not difficult to check that if $\GG$ is $r$-wise $t$-intersecting, 
then so is $\sigma_{i,j}(\GG)$. 
Therefore if $\GG$ is an $r$-wise $t$-intersecting family, 
then there is a shifted $r$-wise $t$-intersecting family $\GG'$ with 
$\mu_p(\GG')=\mu_p(\GG)$.

For $G,H\subset[n]$ we say that $G$ shifts to $H$, denoted by $G\leadsto H$,
if $|G|\leq|H|$ and the $i$th smallest element of $G$ is greater than or 
equal to that of $H$ for each $i\leq|G|$. In this case the walk of $G$ is 
in the lower right area with respect to the walk of $H$.
Note that if $G\leadsto H$ but $|G|\neq|H|$, then 
$H$ cannot be obtained from $G$ by shifting operations only.

We say that $\GG$ is inclusion maximal if $G\in\GG$ and $G\subset H$ imply 
$H\in\GG$. Since we are interested in the maximum measure of 
$r$-wise $t$-intersecting families, we always assume that the families are 
inclusion maximal.
If $\GG$ is shifted and inclusion maximal, then 
$G\in\GG$ and $G\leadsto H$ imply $H\in\GG$.

For a family $\HH$, if there exists $H_0\in\HH$ such that
$H\leadsto H_0$ for all $H\in\HH$, then we say that $H_0$ is the 
shift-end in $\HH$. It depends on the choice of $\HH$ whether the shift-end
in $\HH$ exists or not.

For integers $i\leq j$ let $[i,j]:=\{i,i+1,\ldots,j\}$.
For $1\leq a\leq n$ we use $[a,n]_3$ to denote the set
\begin{align}\label{interval3}
 \bigcup_{i=0}^\infty\{a+3i,a+3i+1\}\cap[n], 
\end{align}
e.g., $[4,10]_3=\{4,5,7,8,10\}$.

\begin{lemma}\label{la1}
Let $\GG\subset 2^{[n]}$ be a shifted $3$-wise $t$-intersecting family.
Then every $G\in\GG$ hits the line $y=2x+t$.
\end{lemma}

\begin{proof}
Suppose, to the contrary, that there is $G\in\GG$ which does not hit the line
$y=2x+t$. Let $H:=[t-1]\sqcup[t+1,n]_3$. Then $H$ is the shift-end of the family
consisting of all subsets which do not hit the line. Thus $G\leadsto H$, and
$H\in\GG$ by the shiftedness of $\GG$. Let $H'=[t]\sqcup[t+2,n]_3$ and
$H''=[t+1]\sqcup[t+3,n]_3$. Since $\GG\ni H\leadsto H'\leadsto H''$, we have
$H,H',H''\in\GG$, but $|H\cap H'\cap H''|=t-1$, contradicting the assumption
that $\GG$ is 3-wise $t$-intersecting.
\end{proof}

\begin{lemma}\label{la2}
Let $0<p<\frac23$ and let $s$ be a positive integer.  
Let $\GG\subset 2^{[n]}$. If every $G\in\GG$ hits the line $y=2x+s$ then 
$\mu_p(\GG)\leq\alpha^s$, where
\[
 \alpha:=\frac12\left(\sqrt{\frac{1+3p}{1-p}}-1\right).
\]
\end{lemma}

\begin{proof}
Consider an infinite random walk in the plane starting from the origin, 
each step of which is a random variable, independent of other steps, 
going up (from $(x,y)$ to $(x,y+1)$) with probability $p$ and 
right (from $(x,y)$ to $(x+1,y)$) with probability $q=1-p$.
Let $P(s)$ be the probability that the random walk hits the line $y=2x+s$.
After the first step, the walk is at $(0,1)$ with probability $p$ and
$(1,0)$ with probability $q$. Thus we have $P(s)=pP(s-1)+qP(s+2)$.
The characteristic equation $x=p+qx^3$ has roots $\alpha,\beta,1$, where
$\beta$ is the same as $\alpha$ with the square root negated.
Then we can write $P(s)=A\alpha^s+B\beta^s+C$ for some constants $A,B,C$.
Since $|\beta|>1$, taking the limit $s\to\infty$, we see that $B=0$.
Also $C=0$ follows from $P(s)\to 0$.
To see this, observe that on average, a step reduces $y-2x$
by $2-3p$, which is positive because $p<\frac23$.
Finally, $P(0)=1$, and so $A=1$. Consequently $P(s)=\alpha^s$.

Now we consider $\mu_p(\GG)$. This is precisely the probability that the random
walk hits the line in the first $n$ steps. 
Thus we have $\mu_p(\GG)\leq \alpha^s$.
\end{proof}

\begin{lemma}\label{la3}
Let $\GG\subset 2^{[n]}$ be a shifted $3$-wise $t$-intersecting family.
For all $F,G,H\in\GG$ there exists $i$ such that 
$|F\cap[i]|+|G\cap[i]|+|H\cap[i]|\geq 2i+t$.
\end{lemma}
\begin{proof}
Suppose, to the contrary, that there is a triple of witnesses
$F,G,H\in\GG$ which satisfies the opposite inequality for every $i$. 
Choose the triple so that $|F\cap G\cap H|$ is minimum. 
Let $j$ be the $t$th element of $F\cap G\cap H$. Then, 
$|F\cap G\cap H\cap[j]|=t$ and
\[
|F\cap[j]|+|G\cap[j]|+|H\cap[j]|< 2j+|F\cap G\cap H\cap[j]|. 
\]
This means that there exists $x\in[j-1]$ such that 
\[
 |F\cap\{x\}| +  |G\cap\{x\}| +  |H\cap\{x\}| \leq 1.
\]
Without loss of generality let $x\not\in G\cup H$, and let
$G'=(G\setminus\{j\})\sqcup\{x\}$. Then by the shiftedness $G'\in\GG$, and
the triple $F,G',H$ is also a witness because $|G'\cap[j]|=|G\cap[j]|$.
But $|F\cap G'\cap H|<|F\cap G\cap H|$, which is a contradiction.
\end{proof}

\begin{lemma}\label{la4}
Let $\GG,\HH\subset 2^{[n]}$. If $\GG$ is shifted and inclusion maximal, and
$H\in\HH\setminus\GG$ is the shift-end in $\HH$, then $\GG\cap\HH=\emptyset$.
\end{lemma}

\begin{proof}
Suppose, to the contrary, that $\GG\cap\HH\neq\emptyset$, and let
$H'\in\GG\cap\HH$. Since $H$ is the shift-end in $\HH$ and $H'\in\HH$, 
we have $H'\leadsto H$. Then, since $H'\in\GG$ and 
$\GG$ is shifted and inclusion maximal, we have $H\in\GG$, a contradiction.
\end{proof}

The following result is well-known, see, e.g., Exercise 5.3.5 (b) in \cite{GJ}, 
and \cite{NT} for some extensions.

\begin{lemma}\label{la5}
The number of walks from the origin to $(s,2s+t)$ which do not cross the
line $y=2x+t$, that is, do not touch the line $y=2x+t+1$, is given by
\[
 f(s,t):=\frac{t+1}{3s+t+1}\binom{3s+t+1}s.
\]
\end{lemma}

The last lemma will be used to determine the extremal configurations
in Theorem~\ref{thm1}.

\begin{lemma}\label{la6}
Let $\GG\subset 2^{[n]}$ be a $3$-wise $t$-intersecting family.
If $\sigma_{i,j}(\GG)=\FF_s^t$ for some $1\leq i<j\leq n$,
then $\GG\cong\FF_s^t$.
\end{lemma}
\begin{proof}
 If $i,j\in[t+3s]$ or $i,j\in[n]\setminus[t+3s]$, then $\GG=\FF_s^t$.
So without loss of generality we may assume that $i=t+3s$ and $j=n$.
Define two subfamilies $\GG_1$ and $\GG_2$ of $\GG$ by
\begin{align*}
\GG_1 &=\{G\in\GG:|G|=t+2s,\,i\not\in G,\, j\in G,\, 
(G\cup\{i\})\setminus\{j\}\not\in\GG\} ,\\
\GG_2 &=\{G\in\GG:|G|=t+2s,\,i\in G,\, j\not\in G,\, 
(G\cup\{j\})\setminus\{i\}\not\in\GG\}.
\end{align*}
Since $\sigma_{i,j}(\GG)=\FF_s^t$ we have $|G\cap[t+3s-1]|=t+2s-1$ for all
$G\in\GG_1\sqcup\GG_2$. If $\GG_1=\emptyset$ then 
$\GG=\sigma_{i,j}(\GG)=\FF_s^t$. If $\GG_2=\emptyset$ then 
$\GG=\{G\subset[n]:|G\cap([i-1]\cup\{j\})|\geq t+2s\}\cong\FF_s^t$.
So the remaining case is $\GG_1\neq\emptyset$ and $\GG_2\neq\emptyset$,
and we will show that this case cannot happen.
If $s=0$ then $|G_1\cap G_2|\leq t-1$ for $G_1\in\GG_1$ and $G_2\in\GG_2$,
a contradiction. So we may assume that $s\geq 1$.

Let $V=\binom{[t+3s-1]}s$. For $H\in V$ let $\bar H:=[t+3s-1]\setminus H$,
and so $|\bar H|=t+2s-1$. Then for every $H\in V$ we have
$\bar H\cup\{j\}\in\GG_1$ or $\bar H\cup\{i\}\in\GG_2$ (but not both).
We define a 3-uniform hypergraph $\HH\subset\binom V3$, by letting
$\{H_1,H_2,H_3\}$ be an edge if $|H_1\cup H_2\cup H_3|=3s$, 
or equivalently, $|\bar H_1\cap \bar H_2\cap \bar H_3|=t-1$.
Let $H,H'\in V$ and 
suppose that $G:=\bar H\cup\{j\}\in\GG_1$ and $G':=\bar H'\cup\{i\}\in\GG_2$. 
Choose $H''\in V$ such that $H''\subset [t+3s-1]\setminus(H\cup H')$
arbitrarily, and let $G'':=\bar H''\cup\{k\}\in\GG$, where $k\in\{i,j\}$. 

We can find $F,F'\in V$ such that both $\{H,F,H''\}$ and $\{H',F',H''\}$ 
are in $\HH$. Let $\tilde F,\tilde F'\in\GG$ be such that 
$\tilde F=F\cup\{f\}$ and $\tilde F'=\bar F'\cup\{f'\}$ for some
$f,f'\in\{i,j\}$. 
If $k=i$ then $|G\cap \tilde F\cap G''|=|H\cap F\cap H''|=t-1$, and 
similarly if $k=j$ then $|G'\cap \tilde F'\cap G''|=t-1$. Thus we get a 
contradiction in either case. 
\end{proof}

\section{Proof of Theorem~\ref{thm1}}
Let $\GG\subset 2^{[n]}$ be a 3-wise $t$-intersecting family,
and let $0<p\leq p_0(t)$.
The result is straightforward if $n=t$, and we may assume that
$n\geq t+1$. Indeed we may assume that $n$ is sufficiently large compared
to $t$, say, $n>t+10$, as explained below.
\begin{lemma}
 Suppose that Theorem~\ref{thm1} holds for some $n=n_0>t$. Then 
Theorem~\ref{thm1} also holds for all $n'$ with $t\leq n'\leq n_0$.
\end{lemma}
\begin{proof}
Let $\GG'\subset 2^{[n']}$ be a 3-wise $t$-intersecting family. 
Define $\GG_0\subset 2^{[n_0]}$ by
\[
 \GG_0:=\GG'\sqcup\{G'\sqcup G_0:G'\in\GG',\,G_0\in 2^{[n'+1,n_0]}\}.
\]
Then $\GG_0$ is also 3-wise $t$-intersecting, and 
$\mu_p(\GG')=\mu_p(\GG_0)$. Applying Theorem~\ref{thm1} to $\GG_0$ we
see that $\mu_p(\GG_0)\leq p^t$, and so $\mu_p(\GG')\leq p^t$.
Moreover, it follows from the construction of $\GG_0$ that
if $\GG_0\cong\FF_0^t(n_0)$ or $\FF_1^t(n_0)$, then
$\GG'\cong\FF_0^t(n')$ or $\FF_1^t(n')$ as well.
\end{proof}

First we show the following.
\begin{prop}
If $t\geq 43$ then $\mu_p(\GG)\leq p^t$.
\end{prop}
\begin{proof}
We may assume that $\GG$ is shifted and inclusion maximal. Let
\[
\lambda:=\max\{i:\text{all walks in $\GG$ hit the line $y=2x+i$}\}.
\]
By Lemma~\ref{la1} we have $\lambda\geq t$. If $\lambda\geq t+1$ then
by Lemma~\ref{la2} we have
\[
 \mu_p(\GG)\leq \alpha^{t+1}<p^t
\]
for $t\geq 9$ (see Claim~\ref{app2} in Appendix for detailed computation).

Thus we may assume that $\lambda=t$. Then all walks in $\GG$ hit the line
$L:y=2x+t$, and some walks do not hit $L':y=2x+t+1$. 
So we divide $\GG$ into three subfamilies:
\[
 \GG=\dot\GG\sqcup\ddot\GG\sqcup\tilde\GG,
\]
where walks in $\tilde\GG$ hit the line $L'$, walks in $\dot\GG$ hit $L$ 
only once, and walks in $\ddot\GG$ hit $L$ at least twice.

Let $\WW$ be the set of all $n$-step walks starting from the origin, 
that is, $\WW=2^{[n]}$, and
we divide $\WW=\dot\WW\sqcup\ddot\WW\sqcup\tilde\WW$ in the same manner.
Since all walks in $\tilde W$ hit the line $L'$, we can apply 
Lemma~\ref{la2} and get
\[
 \mu_p(\tilde W)\leq \alpha^{t+1}.
\]
\begin{claim}
 We have $\mu_p(\ddot\WW)\leq\alpha^{t+2}$.
\end{claim}
\begin{proof}
Let $\WW'\subset 2^{[n]}$ be the set of walks hitting the line $L'':y=2x+t+2$.
It suffices to show that there is an injection from $\ddot\WW$ to $\WW'$.
Let $W\in\ddot\WW$.
Let $P_1$ and $P_2$ be the first and the second points where $W$ hits $L$. 
Reflect the part of $W$ between $P_1$ and $P_2$ over $L$
(and keep the remaining part the same). 
The resulting walk $W'$ hits $L''$, that is, $W'\in\WW'$.
\end{proof}

\begin{claim}\label{claim2}
If $\dot\GG=\emptyset$ then $\mu_p(\GG)<p^t$ for $t\geq 15$.
\end{claim}
\begin{proof}
We have
\[
 \mu_p(\GG)=\mu_p(\ddot\GG\sqcup\tilde\GG)\leq
\mu_p(\ddot\WW)+\mu_p(\tilde\WW)\leq\alpha^{t+2}+\alpha^{t+1}<p^t
\] 
for $t\geq 15$ (see Claim~\ref{app3}).
\end{proof}

Thus we may assume that $\dot\GG\neq\emptyset$.

\begin{claim}
There exists a unique $s$ such that $\dot\GG\sqcup\ddot\GG\subset\FF_s^t$.
\end{claim}

\begin{proof}
Let $G\in\dot\GG$, and suppose that $G$ hits the line $y=2x+t$ only at 
$(s,2s+t)$. Then $G\in\FF_s^t$, and if $s'\neq s$ then $G\not\in\FF_{s'}^t$.
In other words $|G\cap[3s+t]|=2s+t$ and $|G\cap[3s'+t]|\neq 2s'+t$ 
if $s'\neq s$.
(Recall that all walks in $\dot\GG$ hit the line only once, and so $s$
is determined uniquely.)

Choose $G'\in\dot\GG\sqcup\ddot\GG$ arbitrarily. We first show that 
$|G'\cap[j]|\leq\frac{2j+t}3$ for all $j$. Suppose the contrary and 
$|G'\cap[j]|>\frac{2j+t}3$ for some $j$. Then $G'$ hits the point $(x_0,y_0)$ 
where $j=x_0+y_0$ with $y_0>\frac{2j+t}3$. Since
$x_0=j-y_0<j-\frac{2j+t}3=\frac{j-t}3$ and the point 
$(\frac{j-t}3,\frac{2j+t}3)$ is on the line $y=2x+t$, it follows that the point
$(x_0,y_0)$ is strictly above the line, which contradicts the assumption
$G'\in\dot\GG\sqcup\ddot\GG$.

On the other hand, by Lemma~\ref{la3}, there is an $i$ such that 
$|G\cap[i]|+|G'\cap[i]|+|G''\cap[i]|\geq 2i+t$ for all 
$G',G''\in\dot\GG\sqcup\ddot\GG$ (and the $G$ that we have already chosen
in the first paragraph). Then using the fact in the previous
paragraph we have $|G\cap[i]|=|G'\cap[i]|=|G''\cap[i]|=\frac{2i+t}3$.
By the assumption for $G$ we necessarily have $i=3s+t$.
This means $G,G',G''\in\FF_s^t$ and $\dot\GG\sqcup\ddot\GG\subset\FF_s^t$.
\end{proof}  

\begin{claim}
If $t\geq 2$, $s\geq 0$, and $p\leq p_0$, then 
$\mu_p(\FF_s^t)\geq\mu_p(\FF_{s+1}^t)$.
\end{claim}
\begin{proof}
Let $\AA_j=\binom{[t+3s]}{t+2s+j}$, $B=\{t+3s+i:1\leq i\leq 3\}$, 
$C=[t+3s+4,n]$. Then we have
\begin{align*}
\FF_s^t\setminus\FF_{s+1}^t &=
\{F\cup G:F\in\AA_0\cup\AA_1\cup\{A\cup\{b\}:A\in\AA_0,\,b\in B\},\, G\subset C\},\\
\FF_{s+1}^t\setminus\FF_{s}^t &=
\{F\cup B\cup G:F\in\AA_{-1},\,G\subset C\},
\end{align*}
and 
\begin{align*}
\mu_p(\FF_s^t\setminus\FF_{s+1}^t) &=
\tbinom{t+3s}{t+2s}(p^{t+2s}q^{s+3}+3p^{t+2s+1}q^{s+2})
+\tbinom{t+3s}{t+2s+1}p^{t+2s+1}q^{s+2}\\
&\geq\tbinom{t+3s}{t+2s}(p^{t+2s}q^{s+3}+3p^{t+2s+1}q^{s+2}),\\
\mu_p(\FF_{s+1}^t\setminus\FF_{s}^t) &=
\tbinom{t+3s}{t+2s-1}p^{t+2s+2}q^{s+1}.
\end{align*}
Thus $\mu_p(\FF_s^t\setminus\FF_{s+1}^t)\geq\mu_p(\FF_{s+1}^t\setminus\FF_{s}^t)$
if $(s+1)(q^2+3pq)\geq(t+2s)p^2$, that is,
\begin{align}\label{p0max}
 p\leq 2\left(\sqrt{\frac{4t+17s+9}{s+1}}-1\right)^{-1}. 
\end{align}
If $t=2$ then the RHS is $p_0(2)$.
If $t\geq 3$ then the RHS is increasing in $s$, and if moreover $s=0$
then the RHS coincides with $p_0(t)=2/(\sqrt{4t+9}-1)$.
In particular, if $t\geq 2$, $s\geq 0$, and $p\leq p_0$, then we have
$\mu_p(\FF_s^t)\geq\mu_p(\FF_{s+1}^t)$.
\end{proof}

Note that if $t=1$ then the RHS of \eqref{p0max} is decreasing in $s$.
This is why Theorem~\ref{thm1} does not hold for $t=1$.

\begin{claim}
If $s\geq 2$ then $\mu_p(\GG)<p^t$ for $t\geq 43$.
\end{claim}
\begin{proof}
If $s\geq 2$ and $p\leq p_0$ then we have
$\mu_p(\FF_s^t\setminus\FF_{s+1}^t)\geq\mu_p(\FF_{s+1}^t\setminus\FF_{s}^t)$,
and so $\mu_p(\FF_2^t)\geq \mu_p(\FF_3^t)\geq \cdots$.
Note also that $\mu_p(\FF_2^t)=\sum_{i=0}^2\binom{t+6}ip^{t+6-i}q^i$. 
Thus we have
\begin{align*}
 \mu_p(\GG) &= \mu_p(\dot\GG\sqcup\ddot\GG)+\mu_p(\tilde\GG)\\
&\leq\mu_p(\FF_s^t)+\mu_p(\tilde\WW)\\
&\leq\sum_{i=0}^2\binom{t+6}ip^{t+6-i}q^i+\alpha^{t+1}\\
&<p^t
\end{align*}
for $t\geq 43$ (see Claim~\ref{app4}).
\end{proof}

Thus the remaining cases are $s=0$ and $s=1$.
First we deal with the case $s=0$.
Recall the notation $[a,n]_3$ from \eqref{interval3}.

\begin{claim}\label{claim s=0}
If $s=0$ then $\mu_p(\GG)\leq p^t$ for $t\geq 10$.
\end{claim}
\begin{proof}
In this case we may assume that $p\leq p_0(t)\leq\frac13$ because $p_0(t)$ is 
decreasing in $t$ and $p_0(10)=\frac13$.
If $[t]\in\GG$ then $\GG=\FF_0^t$ because $\GG$ is inclusion maximal, 
and we are done. So we may assume that $[t]\not\in\GG$ and we will show 
that $\mu_p(\GG)<p^t$.
For $1\leq i\leq n-(t+1)$ define $W_i\in\dot\WW$ by
\[
W_i:=[t]\sqcup\{t+i+1\}\sqcup[t+i+3,n]_3.
\]
Then $W_1\in\dot\GG$, indeed $W_1$ is the shift-end in $\dot\GG$. 
Thus we can define $I:=\max\{i:W_i\in\dot\GG\}$. 
Since $[t]=W_{n-t}\not\in\GG$, it follows that $I<n-t$.
Let
\[
 \HH:=\{H\in\dot\WW:[t]\subset H\text{ and }H\leadsto W_{I+1}\}\subset\FF_0^t.
\]
Then $W_{I+1}\in\HH\setminus\GG$ is the shift-end in $\HH$, and it follows
from Lemma~\ref{la4} that $\GG\cap\HH=\emptyset$. Thus 
$\FF^t_0\setminus\GG\supset\HH$ and we have
\[
 \mu_p(\FF_0^t\setminus\GG)\geq\mu_p(\HH)\geq p^tq^{I+1}(1-\alpha^2).
\]
Indeed every walk in $\HH$ hits all of $(0,t)$ and $(I+1,t)$; 
and after $t+I+1$ steps it does not hit the line $y=2x+t-2I$, 
see Figure~\ref{fig1}.
(By translation the latter part of the walk can be treated as the walk from
the origin not hitting the line $y=2x+2$.)
\begin{figure}[h]
\includegraphics{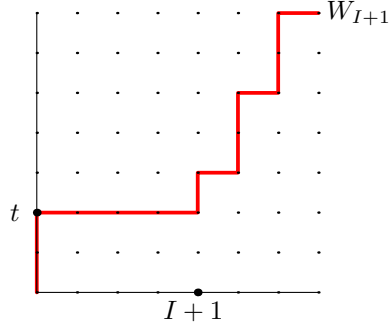}  
\caption{The walk $W_{I+1}$ ($s=0$)}\label{fig1}
\end{figure}

On the other hand we can bound $\mu_p(\GG\setminus\FF_0^t)$ using the 
fact that $W_I\in\GG$. Let 
\[
W':=[t]\sqcup(\{t+I,t+I+2\}\cap[n])\sqcup[t+I+4,n]_3. 
\]
Since $W_I\leadsto W'$ and $\GG$ is shifted, we have $W'\in\dot\GG$. Let 
\[
E:=[t-1] \sqcup([t+1,t+I+3]\cap[n])\sqcup[t+I+5,n]_3,
\]
see Figure~\ref{fig2}.
\begin{figure}
\includegraphics{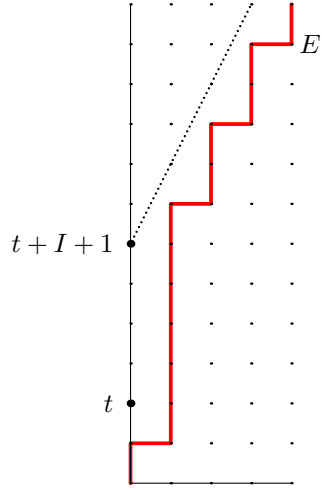}  
\caption{The walk $E$ ($s=0$)}\label{fig2}
\end{figure}
Since $|W_{I}\cap W'\cap E|=t-1$, we have $E\not\in\GG$.
Thus all walks in $\GG\setminus\FF_0^t$ must hit the line $y=2x+t+I+1$, 
and so
$\mu_p(\GG\setminus\FF_0^t)\leq\alpha^{t+I+1}$. (Note that $t+I+1\leq n$.)

We claim that
\begin{align}\label{case:s=0}
\alpha^{t+I+1}<p^tq^{I+1}(1-\alpha^2), 
\end{align}
which implies $\mu_p(\GG\setminus\FF_0^t)<\mu_p(\FF_0^t\setminus\GG)$
and so $\mu_p(\GG)<\mu_p(\FF_0^t)=p^t$.
The inequality \eqref{case:s=0} is equivalent to
$(\alpha/q)^I\alpha^{t+1}<p^tq(1-\alpha^2)$. Since $\alpha/q$ is increasing
in $p$, and at $p=\frac13$ we have $\alpha/q=\frac34(\sqrt3-1)<1$, 
it suffices to show the inequality at $I=1$. 
Thus we need to show that $1<(p/\alpha)^t(q/\alpha)^2(1-\alpha^2)$.
By Claim~\ref{app2} we have $(p/\alpha)^t>\alpha$ and it suffices to show 
that $h(p):=\alpha (q/\alpha)^2(1-\alpha^2)>1$. Indeed $h(p)$ is decreasing
in $p$, and $h(p)>1$ for $p\leq \frac13$. 
\end{proof}

Next we deal with the case $s=1$.

\begin{claim}\label{claim s=1}
If $s=1$ then $\mu_p(\GG)\leq p^t$ for $t\geq 11$.
\end{claim}

\begin{proof}
If $[t+3]\setminus\{t\}\in\GG$ then, by the shiftedness of $\GG$, we have
$G_i:=[t+3]\setminus\{i\}\in\GG$ for all $t\leq i\leq t+3$, 
and $\GG\subset\FF_1^t$. Indeed if $G\in\GG$ satisfies 
$|G\cap[t+3]|<t+2$, say, $[t+3]\setminus G\supset\{i,j\}$,
then we can choose $\{k,l\}\subset[t,t+3]\setminus\{i,j\}$ and
$|G\cap G_k\cap G_l|<t$, a contradiction. 
Moreover, $\GG=\FF_1^t$ follows from the inclusion maximality of $\GG$.

Now assume $[t+3]\setminus\{t\}\not\in\GG$, and we will show that 
$\mu_p(\GG)<\mu_p(\FF_1^t)\leq p^t$. For $1\leq i\leq n-(t+4)$ let 
\[
W_i:=([t+3]\setminus\{t\})\sqcup\{t+i+4\}\sqcup[t+i+6,n]_3\in\dot\WW.
\]
Then $W_1\in\dot\GG$, indeed $W_1$ is the shift-end in $\dot\GG$
because every walk in $\dot\GG$ must hit the line $y=2x+t$ only at $(1,t+2)$.
Thus we can define $I:=\max\{i:W_i\in\dot\GG\}$. 
Since $[t+3]\setminus\{t\}=W_{n-t-3}\not\in\GG$, it follows $I<n-t-3$.
Let
\[
\HH:=\{H\in\dot\WW:|H\cap[t]|=t-1,\,[t+1,t+3]\subset H
\text{ and }H\leadsto W_{I+1}\}\subset\FF_1^t.
\]
Then $W_{I+1}\in\HH\setminus\GG$ is the shift-end in $\HH$, and it follows
from Lemma~\ref{la4} that $\GG\cap\HH=\emptyset$. Thus we have
\[
 \mu_p(\FF_1^t\setminus\GG)\geq\mu_p(\HH)\geq tp^{t+2}q^{I+2}(1-\alpha^2).
\]
Indeed every walk in $\HH$ hits all of $(1,t-1)$, $(1,t+2)$, and 
$(I+2,t+2)$; and after $t+I+4$ steps it does not hit the line $y=2x+t-2I$,
see Figure~\ref{fig3}.
(There are $t$ ways from the origin to $(1,t-1)$ and the first $t+I+4$ steps
to $(I+2,t+2)$ contribute $tp^{t+2}q^{I+2}$.)
\begin{figure}
\includegraphics{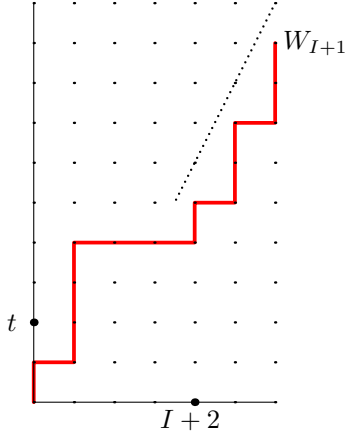}  
\caption{The walk $W_{I+1}$ ($s=1$)}\label{fig3}
\end{figure}

On the other hand we can bound $\mu_p(\GG\setminus\FF_1^t)$ using $W_I\in\GG$.
Let 
\[
W':=([t+3]\setminus\{t+1\})\sqcup(\{t+I+3,t+I+5\}\cap[n])\sqcup[t+I+7,n]_3. 
\]
Since $W_I\leadsto W'$ and $\GG$ is shifted, we have $W'\in\dot\GG$. Let 
\[
E:=[t+1] \sqcup([t+4,t+I+6]\cap[n])\sqcup[t+I+8,n]_3,
\]
see Figure~\ref{fig4}.
\begin{figure}
\includegraphics{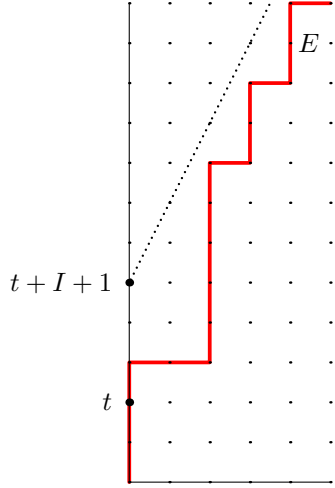}  
\caption{The walk $E$ ($s=1$)}\label{fig4}
\end{figure}
Since $|W_{I}\cap W'\cap E|=t-1$ we have $E\not\in\GG$.
Note that if $G\in\GG\setminus\FF_1^t$ then $|G\cap[t+3]|\leq t+1$.
Thus all walks in $\GG\setminus\FF_1^t$ must hit the line $y=2x+t+I+1$,
and so 
$\mu_p(\GG\setminus\FF_1^t)\leq\alpha^{t+I+1}$. (Note that $t+I+1<n-2$.)

We claim that
\begin{align}\label{case:s=1}
\alpha^{t+I+1}<tp^{t+2}q^{I+2}(1-\alpha^2), 
\end{align}
which implies $\mu_p(\GG\setminus\FF_1^t)<\mu_p(\FF_1^t\setminus\GG)$
and so $\mu_p(\GG)<\mu_p(\FF_1^t)\leq p^t$.
The inequality \eqref{case:s=1} is equivalent to
$(\alpha/q)^{I+1}\alpha^{t}<tp^{t+2}q(1-\alpha^2)$. 
Since $\alpha/q$ is increasing in $p$ and $\alpha/q<1$ for $p\leq p_0(11)$, 
it suffices to show the inequality at $I=1$, that is, 
$t(p/\alpha)^{t+2}q^3(1-\alpha^2)>1$. 
One can verify this inequality by direct computation for $11\leq t\leq 13$.
Let $t\geq 14$. Since $p/\alpha > 1-p^2$ and $q^3(1-\alpha^2)>1-3p+2p^2$, 
we need to show $g(p):=t(1-p^2)^{t+2}(1-3p+2p^2)>1$. 
The LHS is decreasing in $p$ for $p\leq p_0(14)$, 
and so it suffices to show the inequality at $p=p_0(t)$. 
Indeed, $g(p_0(t))$ is increasing in $t$, and $g(p_0(14))>1$, as needed.
\end{proof}

This completes the proof of the proposition.
\end{proof}

To complete the proof of the inequality of Theorem~\ref{thm1} 
we only need to improve the estimation for the case $s\geq 2$.

\begin{claim}\label{claim7}
If $s=2$ then $\mu_p(\GG)<p^t$ for $t\geq 8$.
\end{claim}
\begin{proof}
Let $W\in\dot\GG$ be the shift-end in $\dot\GG$:
\[
W:=([t+8]\setminus\{t,t+3,t+7\})\sqcup[t+10,n]_3, 
\]
see Figure~\ref{fig5}.
\begin{figure}
\includegraphics{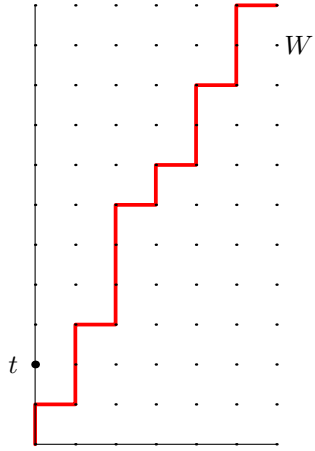}  
\caption{The walk $W$ ($s=2$)}\label{fig5}
\end{figure}
Define $W'\in\GG$ by
\[
 W':=([t+9]\setminus\{t+1,t+4,t+8\})\sqcup[t+11,n]_3.
\]
Then $W\leadsto W'$.
Finally let
\[
 E:=([t+10]\setminus\{t+2,t+5,t+6\})\sqcup[t+12,n]_3,
\]
see Figure~\ref{fig6}.
\begin{figure}
\includegraphics{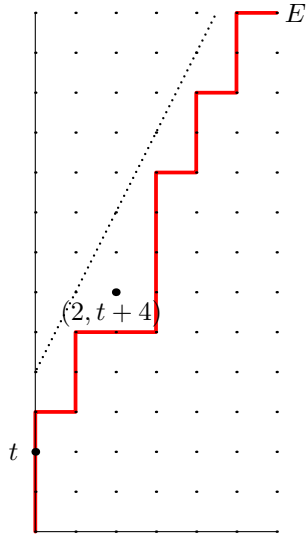}  
\caption{The walk $E$ ($s=2$)}\label{fig6}
\end{figure}

Then $|W\cap W'\cap E|=t-1$ and $E\not\in\GG$.
Thus every walk in $\GG$ hits the line $y=2x+t+2$, or hits $(2,t+4)$.
Note that a walk hitting $(2,t+4)$ without hitting the line must hit
either $(2,t+1)$ or all of $(1,t+1)$, $(1,t+2)$, and $(2,t+3)$.
Therefore we have
\begin{align}\label{eq:claim7}
 \mu_p(\GG) \leq \alpha^{t+2} + \left(\binom{t+3}2+2(t+2)\right)p^{t+4}q^2 <p^t 
\end{align}
for $t\geq 8$ (see Claim~\ref{app5}).
\end{proof}

\begin{claim}\label{claim8}
If $s\geq 3$ then $\mu_p(\GG)<p^t$ for $t\geq 15$. 
\end{claim}

\begin{proof}
 Since $\dot\GG\sqcup\ddot\GG\subset\FF_s^t$ and 
$\GG\subset\FF_s^t\cup\tilde\WW=\tilde\WW\sqcup(\FF_s^t\setminus\tilde\WW)$,
we have
\[
 \mu_p(\GG)\leq\mu_p(\tilde\WW)+\mu_p(\FF_s^t\setminus\tilde\WW).
\]
By Lemma~\ref{la5} it follows
\[
\mu_p(\FF_s^t\setminus\tilde\WW)\leq f(s,t)p^{t+2s}q^s=:g(s,t).
\]
We have
\[
 \frac{f(s,t)}{f(s+1,t)}=\frac{(s+1)(2s+t+3)(2s+t+2)}
{(3s+t+3)(3s+t+2)(3s+t+1)}
>\frac4{t+12}\left(\frac23\right)^2,
\]
and 
\[
\frac{g(s,t)}{g(s+1,t)}=\frac{f(s,t)}{f(s+1,t)p^2 q}>\frac{16}{9(t+12)p^2q}.
\]
We claim that the RHS is $>1$, which means that $g(s,t)$ is decreasing in 
$s$. To this end we show that $\frac{16}9>(t+12)p^2q=:h(p,t)$. 
Since $p^2q$ is increasing in $p$ (for $p<\frac23$), we have 
$h(p,t)\leq h(p_0(t),t)=:\tilde h(t)$. Then $\tilde h(8)<\frac{16}9$ and 
$\tilde h(t)$ is decreasing in $t$ for $t\geq 8$, as needed. 
Indeed 
\[
 \frac d{dt}\tilde h(t)=
-\frac{48\left(3 \sqrt{4t+9}-13\right)}
{\sqrt{4t+9} \left(\sqrt{4t+9}-1\right)^4},
\]
and the RHS is negative if $t\geq 3$.
Thus we have
\[
\mu_p(\FF_s^t\setminus\tilde\WW)\leq g(s,t)\leq g(3,t)\leq
\frac{t+1}{t+10}\binom{t+10}3p^{t+6}q^3=\frac16(t+1)(t+8)(t+9)p^{t+6}q^3.
\]
Consequently we have
\begin{align}\label{eq:claim8}
 \mu_p(\GG)\leq \alpha^{t+1}+\frac16(t+1)(t+8)(t+9)p^{t+6}q^3<p^t
\end{align}
for $t\geq 15$ (see Claim~\ref{app6}).
\end{proof}

We have proved the inequality $\mu_p(\GG)\leq p^t$ with equality holding 
only if $\GG=\FF_0^t$ (Claim~\ref{claim s=0}) or $\GG=\FF_1^t$
(Claim~\ref{claim s=1}). Recall from \eqref{def:p0} that 
$\mu_p(\FF_0^t)=\mu_p(\FF_1^t)$ if and only if $p=p_0$.
Recall also that we have assumed that $\GG$ is shifted. Note that dropping
this assumption does not affect the inequality. 
Now let $\GG$ be a 3-wise $t$-intersecting
family which is not necessarily shifted. We already know that 
$\mu_p(\GG)\leq p^t$. Now suppose that $\mu_p(\GG)=p^t$. 
Starting from $\GG$ we get a shifted family $\GG'$
by applying shifting operations repeatedly. Then $\GG'=\FF_0^t$, or
$\GG'=\FF_1^t$ and $p=p_0$. By Lemma~\ref{la6} we have $\GG\cong\FF_0^t$ or 
$\GG\cong\FF_1^t$. This completes the proof of Theorem~\ref{thm1}.
\qed

\section{Proof of Theorem~\ref{thm2}}\label{sec4}
Let $t\geq 15$ and $0<p\leq p_0$ be fixed.
Let $\GG\subset 2^{[n]}$ be a shifted $t$-intersecting family.
We will basically repeat the proof of Theorem~\ref{thm1} more carefully.

In the proof of Claim~\ref{claim2} we have that if $\dot\GG=\emptyset$ then 
$\mu_p(\GG)\leq\alpha^{t+2}+\alpha^{t+1}<p^t$. So letting 
$\epsilon_1=\epsilon_1(p,t):=p^t-(\alpha^{t+2}+\alpha^{t+1})$ we have 
$\mu_p(\GG)<p^t-\epsilon$ for all $0<\epsilon<\epsilon_1$.

Now we introduce $s$ as in the proof of Theorem~\ref{thm1}. Let $s\geq 2$.
Then, by \eqref{eq:claim7} in Claim~\ref{claim7} and
\eqref{eq:claim8} in Claim~\ref{claim8}, we can choose 
$\epsilon_2$ so that $\mu_p(\GG)<p^t-\epsilon$ for all 
$0<\epsilon<\epsilon_2$. 

Let $\epsilon_0=\epsilon_0(p,t):=\min\{\epsilon_1,\epsilon_2\}$.
Then, except for the cases $s=0$ and $s=1$, it follows that
$\mu_p(\GG)<p^t-\epsilon$ for all $0<\epsilon<\epsilon_0$. 
The remaining cases ($s=0$ and $s=1$) are essential for the stability
as we will see.

First let $s=0$. If $[t]\in\GG$ then $\GG\subset\FF_0^t$ and we are done.
\begin{claim}\label{claim9}
Let $s=0$  and $[t]\not\in\GG$. If $\mu_p(\GG)=p^t-\epsilon$ for some
$0<\epsilon<\epsilon_0$, then $\mu_p(\FF_0^t\triangle\GG)<C_1\,\epsilon$,
where $C_1$ depends only on $p$ and $t$.
\end{claim}

\begin{proof}
 We follow the proof of Claim~\ref{claim s=0}. Let 
$a=\mu_p(\GG\setminus\FF_0^t)$ and $b=\mu_p(\FF_0^t\setminus\GG)$.
We have shown that
\[
 b\geq p^t q^{I+1}(1-\alpha^2)>\alpha^{t+I+1}\geq a,
\]
and 
\[
 \epsilon=\mu_p(\FF_0^t)-\mu_p(\GG)= b-a>
p^t q^{I+1}(1-\alpha^2)-\alpha^{t+I+1}>0.
\]
Indeed,
\begin{align}\label{ineq1:a/b}
 \frac ab\leq\frac{\alpha^{t+I+1}}{p^t q^{I+1}(1-\alpha^2)}
=\left(\frac{\alpha}q\right)^I\frac{\alpha^{t+1}}{p^t q(1-\alpha^2)}
\leq\left(\frac{\alpha}q\right)\frac{\alpha^{t+1}}{p^t q(1-\alpha^2)}<1,
\end{align}
and $\frac ab<1-\delta_1$ for some $\delta_1=\delta_1(p,t)>0$.
It then follows that $\epsilon\geq b-a>\delta_1 b$ and
\[
 \mu_p(\FF_0^t\triangle\GG)=a+b<2b\leq\frac 2{\delta_1}\epsilon.
\]
This means that if $\mu_p(\GG)=p^t-\epsilon$ for some $\epsilon<\epsilon_0$
then $\mu_p(\FF_0^t\triangle\GG)<C_1\,\epsilon$, where
$C_1=\frac2{\delta_1}$ depends only on $p$ and $t$.
\end{proof}

Next let $s=1$. If $[t+3]\setminus\{t\}\in\GG$ then 
$\GG\subset\FF_1^t$ and we are done.
\begin{claim}\label{claim10}
Let $s=1$  and $[t+3]\setminus\{t\}\not\in\GG$. 
If $\mu_p(\GG)=p^t-\epsilon$ for some $0<\epsilon<\epsilon_0$, then 
$\mu_p(\FF_1^t\triangle\GG)<C_2\,\epsilon$, where $C_2$ depends only 
on $p$ and $t$.
\end{claim}

\begin{proof}
We follow the proof of Claim~\ref{claim s=1}. 
Let 
$a=\mu_p(\GG\setminus\FF_1^t)$ and $b=\mu_p(\FF_1^t\setminus\GG)$.
We have shown that
\[
 b\geq tp^{t+2} q^{I+2}(1-\alpha^2)>\alpha^{t+I+1}\geq a,
\]
and 
\begin{align}\label{eq:claim10}
 \epsilon=p^t-\mu_p(\GG)\geq \mu_p(\FF_1^t)-\mu_p(\GG)\geq b-a>
tp^{t+2} q^{I+2}(1-\alpha^2)-\alpha^{t+I+1}>0.
\end{align}
Indeed,
\begin{align}\label{ineq2:a/b}
 \frac ab\leq\frac{\alpha^{t+I+1}}{tp^{t+2} q^{I+2}(1-\alpha^2)}
=\left(\frac{\alpha}q\right)^{I+1}\frac{\alpha^{t}}{tp^{t+2} q(1-\alpha^2)}
\leq\left(\frac{\alpha}q\right)^{2}\frac{\alpha^{t}}{tp^{t+2} q(1-\alpha^2)}
<1,
\end{align}
and $\frac ab<1-\delta_2$ for some $\delta_2=\delta_2(p,t)>0$.
It then follows that $\epsilon\geq b-a>\delta_2 b$ and
\[
 \mu_p(\FF_1^t\triangle\GG)=a+b<2b\leq\frac 2{\delta_2}\epsilon.
\]
This means that if $\mu_p(\GG)=p^t-\epsilon$ for some $\epsilon<\epsilon_0$
then $\mu_p(\FF_1^t\triangle\GG)<C_2\,\epsilon$, where
$C_2=\frac2{\delta_2}$ depends only on $p$ and $t$.
\end{proof}

Let $C=\max\{C_1,C_2\}$.
By Claim~\ref{claim9} and \ref{claim10} it follows that if
$\mu_p(\GG)=p^t-\epsilon$ for some $\epsilon<\epsilon_0$ then 
$\mu_p(\FF\triangle\GG)<C\epsilon$, where $\FF\in\{\FF_0^t,\FF_1^t\}$.
This completes the proof of Theorem~\ref{thm2}. \qed

\medskip
By taking $\epsilon_0$ smaller we obtain another stability result.
To see this let $\epsilon_3=\mu_p(\FF_0^t)-\mu_p(\FF_1^t)$ and let 
\[
 \epsilon_0'=\epsilon_0'(p,t):=
\begin{cases}
\min\{\epsilon_0,\epsilon_3\}&\text{if }p<p_0,\\
\epsilon_0 &\text{if }p=p_0.
\end{cases}
\]
Then $\epsilon_0'>0$ for all $0<p\leq p_0$. (Note that $\epsilon_3=0$ if
$p=p_0$.)

\begin{claim}\label{claim11}
Let $s=1$  and $[t+3]\setminus\{t\}\not\in\GG$. 
If $\mu_p(\GG)=p^t-\epsilon$ for some $0<\epsilon<\epsilon_0'$, then 
$p=p_0$ and $\mu_p(\FF_1^t\triangle\GG)<C_2\,\epsilon$, where 
$C_2$ depends only on $p$ and $t$.
\end{claim}

\begin{proof}
First suppose that $p<p_0$. By \eqref{eq:claim10} we have
\[
p^t-\epsilon=\mu_p(\GG)< \mu_p(\FF_1^t)= p^t-\epsilon_3,
\]
and $\epsilon_3<\epsilon$. This contradicts the assumption that
$\epsilon< \epsilon_0'\leq\epsilon_3$.

Next suppose that $p=p_0$. Then, exactly as in the proof of 
Claim~\ref{claim10}, we have $\mu_p(\FF_1^t\triangle\GG)<C_2\,\epsilon$.
\end{proof}

Using Claim~\ref{claim11} instead of Claim~\ref{claim10} we get the
following.

\begin{thm}\label{thm3}
Let $t\geq 15$ and $0<p\leq p_0$.
Then there exist constants $\epsilon_0'=\epsilon_0'(p,t)>0$ and
$C=C(p,t)>0$ satisfying the following statement: for every 
$0<\epsilon<\epsilon_0'$ and every shifted $3$-wise $t$-intersecting family
$\GG\subset 2^{[n]}$ with $\mu_p(\GG)=p^t-\epsilon$, it follows that
$\mu_p(\FF_0^t\triangle\GG)<C\epsilon$, or $p=p_0$ and
$\mu_p(\FF_1^t\triangle\GG)<C\epsilon$.
\end{thm}

\section{Proof of Theorem~\ref{thm4}}
We modify the proof of Theorem~\ref{thm2}.
Choose $t_0$ so that if $t\geq t_0$ and $p=p_0$ then
\begin{align}\label{ineq:<3/t}
\max\left\{
\left(\frac{\alpha}q\right)^2
\left(\frac{\alpha}p\right)^t
\frac1{1-\alpha^2},\,
\frac1t
\left(\frac{\alpha}p\right)^t
\frac{\alpha^2}{p^2q^3(1-\alpha^2)}
\right\}
< \frac 3t.
\end{align}
This is possible because the LHS is
$\frac et+O((1/t)^{\frac32})$ for $p=p_0$ and $t\to\infty$.
Let $t\geq t_0$ and $p\leq p_0$ be fixed. 
Let $\epsilon=\epsilon_0(p,t)$ be from Theorem~\ref{thm2}.

Let $\GG\subset 2^{[n]}$ be a shifted 3-wise intersecting family with
$\mu_p(\GG)=p^t-\epsilon$, where $0<\epsilon<\epsilon_0$.
We may assume that $\GG$ is inclusion maximal.
If $\GG\subset\FF_0^t$ then $\mu_p(\GG\setminus\FF_0^t)=0$ and we are done.
Thus we may assume that $\GG\not\subset\FF_0^t$, and so
$[t]\not\in\GG$. This corresponds to the case $s=0$ (Claim~\ref{claim9}).
In the same way, we may assume that $\GG\not\subset\FF_1^t$
and $[t+3]\setminus\{t\}\not\in\GG$. This corresponds to the case $s=1$
(Claim~\ref{claim10}).

First suppose that $[t]\not\in\GG$.
Let $a=\mu_p(\GG\setminus\FF_0^t)$ and $b=\mu_p(\FF_0^t\setminus\GG)$.
Note that $\epsilon=b-a$. By \eqref{ineq1:a/b} we have
\[
 \frac ab
\leq\left(\frac{\alpha}q\right)^2
\left(\frac{\alpha}p\right)^t
\frac1{1-\alpha^2}.
\]
The RHS is maximized at $p=p_0$, since $\alpha/q$, $\alpha/p$, and 
$\frac1{1-\alpha^2}$ are all increasing in $p$ for $p\leq p_0$.
Thus by \eqref{ineq:<3/t} we have
$a<\frac3tb=\frac3t(\epsilon+a)$, and $a<\frac{3\epsilon}{t-3}$.

Next suppose that $[t+3]\setminus\{t\}\not\in\GG$.
Let $a=\mu_p(\GG\setminus\FF_1^t)$ and $b=\mu_p(\FF_1^t\setminus\GG)$.
By \eqref{ineq2:a/b} we have
\[
 \frac ab
\leq\frac1t
\left(\frac{\alpha}p\right)^t
\frac{\alpha^2}{p^2q^3(1-\alpha^2)}.
\]
The RHS is maximized at $p=p_0$, since $\alpha/p$, $\frac{\alpha^2}{p^2q^3}$
and $\frac1{1-\alpha^2}$ are all increasing in $p$ for $p\leq p_0$.
Again by \eqref{ineq:<3/t} we have $a<\frac{3\epsilon}{t-3}$.
This completes the proof of Theorem~\ref{thm4}.
\qed

\section*{Acknowledgment}
I thank the referees for their careful reading and many helpful suggestions. 
I also thank one of the referees for pointing out the simple and clean
proof of Lemma~2 which I have adopted.
This research was supported by JSPS KAKENHI Grant No.~18K03399.

\section{Appendix}
Recall that $\alpha=\frac12\left(\sqrt{\frac{1+3p}{1-p}}-1\right)$ and
$p_0=p_0(t)=\frac2{\sqrt{4t+9}-1}$.
By solving $p\leq p_0$ for $t$ we have $t\leq\frac{q(1+2p)}{p^2}=:t_0(p)$. 

\begin{claim}\label{app1}
We have $\alpha<p+p^3$ for $0<p<0.56$.
\end{claim}
\begin{proof}
 We have $\alpha<p+p^3$ if
\[
 \frac{1+3p}{1-p}<(2(p+p^3)+1)^2,
\]
that is, $4p^4 (1-2p+p^2-p^3)>0$, which is true for $0<p<0.56$.
\end{proof}

We define
\[
 \beta(t):=\log p_0(t)+(t+1)p_0(t)^2,
\]
which is decreasing in $t$, indeed,
$\beta'(t)=-\frac{8 t}{\sqrt{4 t+9}(\sqrt{4t+9}-1)^3}<0$.
\begin{claim}\label{app1.5}
Let $k,t_1\geq 1$. 
If $\beta(t_1)<-\log k$ then $k\alpha^{t+1}<p^t$ for all $t\geq t_1$ and
$0<p\leq p_0(t)$.
\end{claim}
\begin{proof}
Since $\alpha/p$ is increasing in $p$, it suffices to 
show $p(\alpha/p)^{t+1}<1/k$ at $p=p_0(t)$. 
Using the previous claim we have $p(\alpha/p)^{t+1}<p(1+p^2)^{t+1}$, and 
it suffices to show that
\[
 \log p_0+(t+1)\log(1+p_0^2)<-\log k.
\]
This follows from $\beta(t)<-\log k$ because $\log(1+p_0^2)<p_0^2$.
Since $\beta(t)$ is decreasing in $t$, we only need $\beta(t_1)<-\log k$, 
which is our assumption.
\end{proof}

\begin{claim}\label{app2}
Let $t\geq 9$ and $0<p\leq p_0$. Then $\alpha^{t+1}<p^t$.
\end{claim}
\begin{proof}
It suffices to show that $p(\alpha/p)^{t+1}<1$ at $p=p_0(t)$.
This can be verified for $9\leq t\leq 12$ by direct computation. 
Let $t\geq 13$. Using Claim~\ref{app1.5} for $k=1$ and $t_1=13$, 
we only need to check that $\beta(13)<0$, which is true.
\end{proof}

\begin{claim}\label{app2.5}
Let $t\geq 20$ and $0<p\leq p_0$. Then $1.29\alpha^{t+1}<p^t$.
\end{claim}
\begin{proof}
This inequality follows form the fact $\beta(20)<-\log 1.29$ and
Claim~\ref{app1.5} for $k=1.29$ and $t_1=20$. 
\end{proof}

\begin{claim}\label{app3}
Let $t\geq 15$ and $0<p\leq p_0$. Then $\alpha^{t+2}+\alpha^{t+1}<p^t$.
\end{claim}
\begin{proof}
It suffices to show that 
$p^2(\alpha/p)^{t+2}+p(\alpha/p)^{t+1}<1$ at $p=p_0(t)$.
We can verify this for $15\leq t\leq 19$ directly.
Let $t\geq 20$ and let $k=1.29$.
Using Claim~\ref{app1} with $p_0=p_0(t)\leq p_0(20)<0.24$ we have
\[
 \alpha+1 < p+p^3+1 \leq 1+p_0(20)+p_0(20)^3 < k.
\]
Then by Claim~\ref{app2.5} we have
$\alpha^{t+2}+\alpha^{t+1}= (\alpha+1)\alpha^{t+1} < (\alpha+1)p^t/k< p^t$.
\end{proof}

\begin{claim}\label{app4}
Let $t\geq 43$ and $0<p\leq p_0$. 
Then $\sum_{i=0}^2\binom{t+2}ip^{t+6-i}q^i+\alpha^{t+1}<p^t$. 
\end{claim}
\begin{proof}
Since $\beta(43)<-\log 2$ we have
$2\alpha^{t+1}<p^t$ for $t\geq 43$ by Claim~\ref{app1.5}.
So it suffices to show that 
$\sum_{i=0}^2\binom{t+2}ip^{t+6-i}q^i<\frac12p^t$. Let 
\[
g(p,t):=p^{-t}\sum_{i=0}^2\binom{t+2}ip^{t+6-i}q^i 
=\frac{1}{2} \,p^4
   \left(q^2 t^2-(p-3)qt+2
\right).
\]
Then
\[
\frac{\partial}{\partial t}\, g(p,t)= 
\frac12\,q p^4(2qt+q+2)>0.
\]
Thus $g(p,t)$ is increasing in $t$, and it follows that
\[
g(p,t)\leq g(p,t_0(p))=\frac{1}{2}
\left(p^3-p^2+p+1\right)
   \left(2 p^3-p^2-p+1\right)=:\tilde g(p),
\]
and $\tilde g(p)<\frac12$ for $0<p\leq 0.8$. 
(Note that $p_0\leq p_0(43)<0.161$.)
\end{proof}

\begin{claim}\label{app5}
Let $t\geq 8$ and $0<p\leq p_0$. Then
$\alpha^{t+2} +\left(\binom{t+3}2+2(t+2)\right)p^{t+4}q^2 <p^t$.
\end{claim}
\begin{proof}
 We show that 
$\alpha^{t+2}/p^t+\frac12(t^2+9t+14)p^4q^2<1$.
This can be verified for $t=8$ by direct computation. Let $t\geq 9$.
By Claims~\ref{app1} and \ref{app2} it follows that
$\alpha^{t+2}/p^t<\alpha<p+p^3$.
Thus it suffices to show that 
$h(p,t):=(p+p^3)+\frac12(t^2+9t+14)p^4q^2<1$.
Since $h(p,t)$ is increasing in $t$, we have 
$h(p,t)\leq h(p,t_0(p))=
\frac{1}{2} \left(5 p^5-4p^4-3 p^3+3 p^2+2p+1\right)=:\tilde h(p)$.
Then $\tilde h(p)<1$ for $p\leq p_0(9)$.
\end{proof}

\begin{claim}\label{app6}
Let $t\geq 15$ and $0<p\leq p_0$. Then
$\alpha^{t+1}+\frac16(t+1)(t+8)(t+9)p^{t+6}q^3<p^t$.
\end{claim}
\begin{proof}
We can verify the inequality for $15\leq t\leq 19$ directly. Let $t\geq 20$.
We show that $\alpha^{t+1}/p^t+\frac16(t+1)(t+8)(t+9)p^{6}q^3<1$.
For the first term of the LHS, by Claim~\ref{app2.5}, 
we have $\alpha^{t+1}/p^t<1/1.29<0.7752$. For the second term,
let $h(p,t)=\frac16(t+1)(t+8)(t+9)p^{6}q^3$.
Then $h(p,t)$ is increasing in $t$ and $h(p,t)\leq h(p,t_0(p))=:\tilde h(p)$,
and $\tilde h(p)=\frac16q^3(1+pq)(1+p+6p^2)(1+p+7p^2)$. 
Since $\tilde h(p_0(20))<0.2244$, and $\tilde h(p)$ is increasing in $p$ 
for $p\leq p_0(20)$, we have $\tilde h(p)<0.2244$.
Finally the result follows from $0.7752+0.2244<1$.
\end{proof}


\begin{thebibliography}{99}
\bibitem{AK1}
R.~Ahlswede, L.H.~Khachatrian.
\newblock The complete intersection theorem for systems of finite sets.
\newblock {\em European J.\ Combin.}, 18:125--136, 1997.

\bibitem{AK2}
R.~Ahlswede, L.H.~Khachatrian.
\newblock A Pushing-pulling method: new proofs of intersection theorems.
\newblock {\em Combinatorica}, 19:1--15, 1999.

\bibitem{AK-p}
R.~Ahlswede, L.H.~Khachatrian.
\newblock The diametric theorem in Hamming spaces--optimal anticodes.
\newblock {\em Adv.\ in Appl.\ Math.}, 20:429--449, 1998.

\bibitem{BE}
C.~{Bey}, K.~Engel.
\newblock Old and new results for the weighted $t$-intersection
problem via AK-methods.
\newblock {\em Numbers, Information and Complexity, Althofer, Ingo,
Eds. et al., Dordrecht}, Kluwer Academic Publishers, 45--74, 2000.

\bibitem{DS}
I.~Dinur, S.~Safra.
\newblock On the Hardness of Approximating Minimum Vertex-Cover.
\newblock {\em Annals of Mathematics}, 162:439-485, 2005.

\bibitem{EKL18}
D.~Ellis, N.~Keller, N.~Lifshitz. 
Stability for the Complete Intersection Theorem, and the Forbidden 
Intersection Problem of Erd\H{o}s and S\'os.
preprint.  
arXiv:\href{https://arxiv.org/abs/1604.06135}{1604.06135}

\bibitem{EKL19}
D.~Ellis, N.~Keller, N.~Lifshitz. 
Stability versions of Erd\H{o}s--Ko--Rado type theorems via isoperimetry. 
J.\ Eur.\ Math.\ Soc.\ (JEMS) 21 (2019), no. 12, 3857--3902.

\bibitem{EKR}
P.~{Erd\H os}, C.~Ko, R.~Rado.
\newblock Intersection theorems for systems of finite sets.
\newblock {\em Quart.\ J.\ Math.\ Oxford (2)}, 12:313--320, 1961.

\bibitem{Filmus}
Y.~Filmus. 
The weighted complete intersection theorem. 
J. Combin. Theory Ser. A 151 (2017), 84--101.

\bibitem{FGL}
Y.~Filmus, K.~Golubev, N.~Lifshitz.
High dimensional Hoffman bound and applications in extremal combinatorics. 
Algebr. Comb. 4 (2021), no. 6, 1005--1026.

\bibitem{F76}
P.~Frankl.
On Sperner families satisfying an additional condition. 
J. Combinatorial Theory Ser. A 20 (1976), no. 1, 1--11.

\bibitem{F87}
P.~Frankl.
The shifting technique in extremal set theory.
 Surveys in combinatorics  (New Cross, 1987) 81--110, 
London Math.\ Soc.\ Lecture Note Ser. 123.

\bibitem{Frankl1991}
P.~Frankl. 
Multiply-intersecting families. 
J. Combin. Theory Ser. B 53 (1991), no. 2, 195--234. 

\bibitem{Frankl2019}
P.~Frankl. 
Some exact results for multiply intersecting families. 
J. Combin. Theory Ser. B 136 (2019), 222--248.

\bibitem{FT2003}
P.~Frankl, N.~Tokushige.
Weighted multiply intersecting families.
Studia Sci.\ Math.\ Hungar.\ 40 (2003) 287--291. 

\bibitem{FLST}
P.~Frankl, S.~J.~Lee, M.~Siggers, N.~Tokushige.
An Erd\H{o}s--Ko--Rado theorem for cross $t$-intersecting families.
J.\ Comb.\ Theory (A), vol 128 (2014) 207--249.

\bibitem{Friedgut}
E.~Friedgut.
On the measure of intersecting families, uniqueness and stability.
Combinatorica 28, 503--528 (2008)

\bibitem{GJ}
I.~P.~Goulden, D.~M.~Jackson. 
Combinatorial Enumeration, Dover, Mineola, NY, 2004 
(reprint of 1983 original).

\bibitem{Gronau}
H.-D.~O.~F.~Gronau. 
On Sperner families in which no $k$ sets have an empty intersection. 
III. Combinatorica 2 (1982), no. 1, 25--36.

\bibitem{KS}
G.~Kindler, S.~Safra.
Noise-Resistant Boolean-Functions are Juntas.
preprint (2002)

\bibitem{LST}
S.~J.~Lee, M.~Siggers, N.~Tokushige.
AK-type stability theorems on cross t-intersecting families. 
European J.\ Combin.\ 82 (2019) 102993, 20 pp.

\bibitem{MV}
D.~Mubayi, J.~Verstra\"ete.
Proof of a conjecture of Erd\H{o}s on triangles in set-systems.
Combinatorica 25 (2005), no. 5, 599--614.

\bibitem{NT}
T.~Nakamigawa, N.~Tokushige.
Counting lattice paths via a new cycle lemma.
SIAM J.\ Discrete Math.\ 26 (2012) 745--754.

\bibitem{Tuvsw}
N.~Tokushige.
\newblock Intersecting families --- uniform versus weighted.
\newblock{\em Ryukyu Math.\ J.}, 18 (2005) 89--103.

\bibitem{T2007}
N.~Tokushige.
The maximum size of 3-wise $t$-intersecting families, 
European J.\ of Combin.\ 28 (2007) 152--166.

\bibitem{T2007b}
N.~Tokushige. 
Multiply-intersecting families revisited. 
J. Combin. Theory Ser. B 97 (2007), no. 6, 929--948. 

\bibitem{T2010}
N.~Tokushige. 
A multiply intersecting Erd\H{o}s--Ko--Rado theorem 
--- the principal case. Discrete Math. 310 (2010), no. 3, 453--460. 

\bibitem{T2021}
N.~Tokushige.
Application of hypergraph Hoffman's bound to intersecting families.
Algebr. Comb. 5 (2022), no. 3, 537--557.

\bibitem{T2022}
N.~Tokushige.
The maximum measure of non-trivial 3-wise intersecting families.
preprint.  
arXiv:\href{https://arxiv.org/abs/2203.17158}{2203.17158}
\end{thebibliography}
\end{document}